\documentclass[12pt]{amsart}
\usepackage{textcase}
\usepackage{amsmath,amsthm,,amssymb}
\usepackage[utf8]{inputenc}
\usepackage[T2A]{fontenc}
\usepackage[russian]{babel}
\usepackage[margin=2cm]{geometry}
\usepackage{secdot}
\usepackage[usenames]{xcolor}
\theoremstyle{remark}
\newtheorem{exa}{\textbf{Пример}}

\def\R{\mathbb{R}}
\def\ul{\underline}

\title[Обучение доказательству будущих учителей математики\ldots]{Обучение доказательству будущих учителей математики на занятиях по математическому анализу}
\author{А. Х. Назиев, И. В. Землякова}
\date{}

\begin{document}

\maketitle

\begin{abstract}
{\flushleft В работе рассматривается построение курса математического анализа педвуза, направленное на формирование у будущих учителей математики умения обнаруживать и решать проблемы, связанные с отысканием доказательств и обучения этому.

{\bf Ключевые слова}: преподавание математики, подготовка учителя, открытие доказательств, кванторы, математический анализ}
\end{abstract}

\section{Чему учить будущего учителя математики?}

На этот вопрос могут быть даны разные ответы. Важно, чтобы он ставился и получал ответ: и в сознании каждого отдельного учителя, и в сознании общества в целом.
Наш ответ на этот вопрос гласит:

\begin{center}\em
Будущего учителя математики нужно учить преподаванию математики~\cite{what}.
\end{center}

Этот ответ может показаться тривиальным, однако он не настолько тривиален, насколько кажется. В нём содержится по меньшей мере три ценных указания:
\begin{itemize}
\item будущего учителя математики {\em нужно учить};
\item его нужно учить {\em преподавать}; и
\item преподавать --- именно {\em математику}.
\end{itemize}
Ни одно из этих положений мы бы не назвали общепризнанным.

В справедливости первого сомневается подавляющее большинство студентов педвузов\footnote{то есть, пединститутов, педуниверситетов, педагогических отделений университетов}. Они уверены: всё, чему их учат в вузе, никогда им не пригодится, а то, что действительно нужно, они освоят, когда начнут работать. И их трудно в этом упрекнуть, потому что даже Российский премьер призвал идти работать в школу людей без педагогического образования. Но и без него немалое число людей, в том числе~--- преподавателей математики педвузов, явно или неявно сомневается в справедливости перечисленных положений. Чтобы убедиться в этом, достаточно заглянуть в существующие учебники основных математических дисциплин педвузов и посетить занятия по этим предметам. Мы увидим, что авторы учебников и преподаватели математических дисциплин с увлечением обучают студентов {\em своим предметам}: алгебре и теории чисел, геометрии, математическому анализу и т.~д. Мы же утверждаем, что все они должны обучать будущего учителя математики одному и тому же --- {\em преподаванию математики}, с той лишь разницей, что каждый из них будет делать это на материале своего курса.

Так же, как школьный учитель математики, преподавая математику, должен понимать, что главное для большинства учеников заключается отнюдь не в самой математике, а в той роли, которую она играет в их образовании, так и преподаватель того или иного раздела математики педвуза, обучая студентов своему предмету, должен понимать, что главное для большинства его студентов заключается не в самом этом предмете, а в той роли, которую он сыграет в становлении их как учителей математики. По этому поводу Б.~В.~Гнеденко в 1989 году писал~\cite{gned}:
\begin{quote}
Мне известны многочисленные случаи, когда с первых лекций преподавание курса математического анализа подчинено единственной цели~--- подготовке изложения в идейном отношении последующих глав математики --- многомерного и функционального анализа. Мне кажется, что основная задача курса пединститутов (а также педагогических отделений университетов) состоит в другом --- в подготовке увлечённого, знающего и умеющего преподавать учителя.
\end{quote}

Сформулируем развёрнуто эту установку Б. В. Гнеденко, восстановив опущенные, но подразумеваемые им слова.

\begin{center}
\em Основная задача всех без исключения, а значит --- и математических, дисциплин педвузов состоит в подготовке увлечённого, знающего и умеющего преподавать математику учителя математики.
\end{center}

Это требует от каждого преподавателя математики педвуза так преподавать свой предмет, чтобы средствами этого предмета готовить учителя математики, то есть обучать преподаванию математики. Способ решения этой задачи указан Д. Пойа в его книге "<Математическое открытие">. В четырнадцатой  главе этой книги её автор пишет:
\begin{quote}
Все курсы, которые я читал учителям математики, были построены так, чтобы они могли служить в какой-то мере и курсами методики. В названии курса обычно указывался только учебный предмет, которому посвящался курс, отводимое же время распределялось между математикой и методикой её преподавания: вероятно, девять десятых всего времени тратилось на предмет и одна десятая --- на методику~\cite[С. 304]{Polya:disc}.
\end{quote}

В применении к курсу математического анализа сказанное принимает следующий вид:

\begin{center}
\em курс математического анализа педвуза должен быть построен так, чтобы он мог служить в какой-то мере и курсом методики преподавания математики.
\end{center}

Разумеется, основное внимание должно быть уделено математическому анализу, но какая-то часть ("<одна десятая">) --- методике, то есть педагогической стороне дела. Конечно же, эта сторона должна тщательно прорабатываться любым преподавателем математического анализа, но не на всех направлениях подготовки уместно посвящать в неё обучающихся. При работе же с будущими учителями математики делать это не только уместно, но и необходимо.

Подробному решению указанной задачи мы намерены посвятить не\-сколько последующих работ. В них мы предполагаем показать, как могут быть реализованы высказанные пожелания, причём не в ущерб математическому анализу, а с пользой для него. Здесь же мы ограничимся первыми примерами. Но сначала обратим внимание на одно обстоятельство, которое могло остаться незамеченным. Мы имеем в виду слова \glqq преподавание математики\grqq. Преподавать будущим учителям математики алгебру, геометрию, математический анализ и т. д. нужно так, чтобы, тем самым, учить {\em преподаванию математики}. А для этого нужно точно знать, что такое математика, что значит преподавать её и зачем это нужно делать. Систему согласованных ответов на перечисленные вопросы мы называем концепцией преподавания математики.

\section{Концепция преподавания математики}

Повторим: будущего учителя математики нужно обучать преподаванию математики. Это требует и от того, кто собирается преподавать математику, и от того, кто будет обучать этому, свободного владения предметом рассмотрения: преподаванием математики. А для этого каждому из них необходимы совершенно определённые ответы на вопросы: Что такое математика? Что значит преподавать математику? Зачем это следует делать?
\begin{quote}
Это будет то, что он\footnote{Учитель. --- А.~Н.} должен постоянно держать в своём уме в качестве
момента, определяющего всю его работу; интересы хорошего преподавания требуют, чтобы преподаватель
знал не только то, \ ч е м у \ он учит, знал не только \ к а к \ учить, но и \ з а ч е м \ он
учит~\cite[С.~9]{joung}.
\end{quote}

Систему согласованных ответов на вопросы: Что такое математика? Что значит преподавать математику? Зачем нужно преподавать математику? --- мы называем концепцией преподавания математики. Разумеется, могут быть разные концепции. На настоящий момент (февраль 2021) нам известна только одна~--- та, что сформулирована и всесторонне обоснована первым автором настоящей работы в его докторской диссертации~\cite{dissB} и многочисленных сопутствующих ей и последовавших за ней публикациях (например, в~\cite{humed}). Подчеркнём, что речь идёт о концепции именно {\em преподавания математики}, не о концепции развития математики, математического образования и т. п. Приведём формулировку этой концепции.

\begin{center}
\fbox{\fbox{
\parbox{12cm}{

\centerline{\underline{\sc Концепция преподавания математики}}

\noindent
 1. Математика~--- это искусство доказательства.

\noindent
 2. Преподавать математику~--- значит систематически
    побуждать учащихся к открытию собственных
    доказательств.

\noindent
 3. Преподавание математики является незаменимым средством
    формирования человека культурного: мыслящего, нравственного,
    свободного.}
}}
\end{center}

По сравнению с нашей первоначальной формулировкой мы добавили одно слово: \glqq искусство\grqq. Дело в том, что математику, в отличие от других областей знания, можно характеризовать с двух точек зрения: с точки зрения объекта изучения (как это и делается в случае остальных областей знания) и с точки зрения метода. С точки зрения объекта изучения математика --- наука: наука о количественных отношениях и пространственных формах действительного мира; с точки зрения метода математика --- искусство: искусство доказательства. Может быть, стоило бы дать более развёрнутую формулировку: математика~--- это наука о количественных отношениях и пространственных формах действительного мира, соединённая с искусством доказательства. Невнимание к этой второй стороне математики мешает понять, почему математика трудна для обучения. --- Потому, что она является не только наукой, но и искусством. Науке научить можно, и этой стороне математики удаётся достаточно успешно обучать. А вот искусству научить невозможно, искусством можно только овладеть посредством самостоятельных усилий, достаточно длительных, достаточно интенсивных. Эта-то сторона дела и вызывает наибольшие трудности при обучении математике, как у обучающихся, так и у обучающих. Во-первых, большинство обучающихся не горит желанием нужные усилия прилагать, поэтому, во-вторых, и обучающие не стремятся побуждать их к этому, зная, какое сопротивление будет ответом на их действия. Гораздо легче диктовать готовые решения и заставлять их выучивать. Но к успеху в обучении математике это привести не может, что и подтверждается из года в год на разного рода экзаменах. Поэтому, как бы ни было трудно, нужно прилагать требуемые усилия, и обучающим, и обучающимся. Именно этой стороне дела и посвящена наша работа в отношении курса математического анализа. И здесь предложенная концепция оказывает существенную и действенную помощь, поскольку совершенно точно указывает, какие именно нужно прилагать усилия --- усилия, направленные на поиск доказательств и побуждение к этому.

Отметим, что подчёркивание двойственного характера математики позволяет примирить отечественную и западную традиции в трактовке математики. В отечественной традиции математика --- наука, в западной --- искусство ("<свободное искусство">, liberal art). Предложенная нами точка зрения выявляет односторонний характер и той, и другой традиции (каждая из них подчеркивает лишь одну из сторон математики) и позволяет преодолеть эту односторонность.

Слово \glqq доказательство\grqq\ очень часто вызывает негативную реакцию и у обучающихся, и у их родителей, и у обучающих. Тому есть две причины. 1) В доказательстве видят лишь одно его назначение, средство подтверждение истинности, нечто вроде прокурорского надзора. Отсюда и обычная ученическая реакция на попытки доказательства: "<Не надо, мы Вам верим">.

2) В сознании большинства людей слово \glqq доказательство\grqq\ невольно ассоциируется со словом \glqq теорема\grqq: им говорят \glqq доказательство\grqq, а они слышат \glqq доказательство {\em теорем}\grqq. Против этого не стоило бы возражать, если бы слово \glqq теорема\grqq\ понималось так, как оно понимается в современной (математической) логике, --- как предложение, для которого найдено доказательство. Но не так оно понимается теми, о ком сейчас речь. Говоря о теоремах, они имеют в виду предложения, называемые так в учебниках и выделенные в них жирным шрифтом. А под доказательством теорем понимают выучивание и пересказывание готовых доказательств. Эта сторона процесса обучения математике похожа на издевательство над здравым смыслом и совершенно справедливо вызывает отторжение у большинства изучающих математику.
\begin{quote}\small

--- Докажи мне теорему Пифагора, --- говорит ученику учитель. При этом он отнюдь не ожидает, что ученик сам докажет названную им теорему, а ожидает услышать от ученика заранее выученное им доказательство.

--- Какую теорему? --- переспрашивает ученик.

--- Теорему Пифагора, --- отвечает учитель.

--- А почему она так называется? --- спрашивает ученик.

--- Потому, что это Пифагор её доказал, --- отвечает учитель.

--- Так, он её доказал? --- спрашивает ученик.

--- Ну, да, доказал, --- отвечает учитель.

--- Зачем же я стану её передоказывать? --- недоумевает ученик.
\end{quote}

Ученик абсолютно прав: зачем передоказывать уже доказанное, когда так много встречается недоказанного. Но именно так странно устроено традиционное преподавание математики: насильственное насаждение передоказывания уже доказанного --- и беззаботность в отношении доказательства предложений, о которых неизвестно, доказаны они или нет. Сплошь и рядом ученики бездоказательно пишут: $\sqrt{\sin^2 x} = \pm\sin x$, $\lg xy = \lg x + \lg y$, $\sqrt[3] x = x^\frac13$ и т.~п., зачастую не вызывая нареканий учителей, но при этом зачем-то "<доказывают">\ теоремы Фалеса, Пифагора, Виета и т.~д. Да, что там ученики: в п. 4 мы приведём пример неверного утверждения, сформулированного без доказательства именитыми преподавателями математики из МГУ в написанном ими пособии по математике \glqq для поступающих в серьёзные вузы\grqq!

\section{Ещё раз о том, что такое функция}
Здесь мы, основываясь на~\cite{what}, продемонстрируем, как установка на доказательство позволяет обнаруживать ошибки в ситуациях, казалось бы, далёких от доказательства.

Согласно определению из школьного учебника~\cite[п.~3]{abma}:
\begin{quote}
числовой функцией с областью определения $D$ называется соответствие, при котором каждому числу $x$ из множества $D$ сопоставляется по некоторому правилу число $y$, зависящее от~$x$.
\end{quote}
Задержим внимание на некоторых частях этого определения.

Начнём с упоминания множества $D$. Рассмотрим обычную квадратичную функцию (т.~е. соответствие, которое каждому действительному числу сопоставляет его квадрат и больше ничему ничего не сопоставляет), а в качестве~$D$ возьмём множество, элементами которого являются: числа 0 и 1, Луна и пирамида Хеопса~--- и только они. Является ли указанное соответствие ($x\leadsto x^2\colon \R\to\R$) числовой функцией с областью определения~$D$? Является, ибо каждому числу $x$ из множества $D$ оно сопоставляет по некоторому правилу число $y$, зависящее от $x$. Таким образом, этот пример показывает, что в приведённой формулировке из учебника недосказано, что $D$ состоит {\em лишь} из действительных чисел. Но не только это.

Является ли рассмотренное соответствие числовой функцией с областью определения $\mathbb R$? Является, ибо каждому числу из множества $\mathbb R$ \ldots{} ну, и так далее. И числовой функцией с областью определения $[0, 1]$~--- тоже является, и числовой функцией с областью определения $(\pi,\ 123]$, и вообще~--- числовой функцией с любой областью определения, содержащейся во множестве всех действительных чисел. Ибо, для того, чтобы именно $D$ было областью определения, мало {\em назвать} множество $D$ областью определения функции, нужно ещё {\em обеспечить} это. А для этого недостаточно, чтобы каждому числу из множества $D$ что-то сопоставлялось, нужно ещё, чтобы больше ничему ничего не сопоставлялось, а об этом в приведённом определении ничего не сказано.

Обратим также внимание на слова <<сопоставляется по некоторому правилу>>. Рассмотрим соответствие, которое каждому действительному числу $x$ сопоставляет число $\sin^2 x + \cos^2 x$, а больше ничему ничего не сопоставляет, и соответствие, которое каждому действительному числу сопоставляет число 1, а больше ничему ничего не сопоставляет. Правила, по которым осуществляется сопоставление в первом и втором случаях, совершенно различны, соответствия же, согласно общепринятым в математике представлениям,~--- одинаковы. Для понятия функции неважно, по какому правилу сопоставляется, важно, {\em что} сопоставляется. Правила важны, но не для {\em понятия функции}, а для {\em задания конкретных функций}, поэтому в определении общего понятия функции о них говорить незачем.

Зато нужно подчеркнуть, что каждому числу из области определения сопоставляется ровно одно число. Для различных $x$ это могут быть разные числа, но каждый раз~--- только одно. Авторы учебника, видимо полагают, что об этом достаточно красноречиво свидетельствует единственное число слова <<число>> (простите за невольный каламбур) в выражении <<число $y$>>. Так или иначе~--- важно отдавать себе отчёт в этой однозначности.

И в заключение скажем о словах <<число $y$, зависящее от $x$>>. Ими авторы хотят подчеркнуть, что <<число $y$>>, хотя и единственно для каждого отдельного~$x$, вовсе не обязано быть одним и тем же для всех $x$, а может меняться от одного $x$ к другому, как бы <<зависеть>> от $x$. Слово <<зависеть>> здесь не совсем удачно, ибо, хотя это $y$ и не обязано быть одним и тем же для различных $x$, вполне может оказаться таковым~--- для постоянных функций,~--- а тогда как раз принято говорить, что оно {\em не} зависит от $x$. Так что и от упоминания об этой <<зависимости>> в общем определении функции лучше воздержаться.

Учитывая сделанные замечания, приходим к следующей (предварительной) формулировке (добавленные нами слова подчёркнуты).
\begin{quote}
Пусть $D$~--- какое-нибудь множество \ul{действительных чисел}. Числовой функцией с областью определения $D$ называется соответствие, при котором каждому числу из множества $D$ сопоставляется ровно по одному числу, \ul{и больше ничему} \ul{ничего не сопоставляется}.
\end{quote}

Теперь заметим, что если каждому числу из множества $D$ сопоставляется ровно по одному числу, а больше ничему ничего не сопоставляется, то каждому вообще действительному числу сопоставляется не более одного числа. При этом $D$ однозначно восстанавливается~--- как множество всех тех и только тех действительных чисел, которым что-то сопоставлено,~--- так что упоминать его в определении функции совершенно незачем. Это наблюдение приводит к следующему, окончательному и, на наш взгляд, наиболее простому и естественному определению.

\begin{quote}
{\em Числовой функцией называется соответствие между действительными числами, при котором каждому действительному числу сопоставляется не более одного действительного числа. Областью определения числовой функции $f$ называется множество $D(f)$ всех тех и только тех действительных чисел, которым в соответствии $f$ сопоставляется хотя бы одно действительное число.}
\end{quote}

\medskip
Слова <<соответствие между действительными числами>> означают, что в соответствии {\sc изначально} участвуют лишь действительные числа. Иначе говоря: 1)~каждый объект, которому что-то {\sc может быть} сопоставлено в этом соответствии, есть действительное число; и 2)~каждый объект, который чему-то {\sc может быть} сопоставлен в этом соответствии, есть действительное число. Это, в частности, делает ненужной (в приведённом определении, но не в заданиях конкретных функций!) существенную добавку <<и больше ничему ничего не сопоставляется>>, упущенную в школьном учебнике.

{\sl Заключительное замечание.} Разумеется, определение, к которому мы пришли, <<хорошо известно>>. Однако почему-то не оно используется в большинстве учебников математики, как в школе, так и в вузе. Наша цель состояла в том, чтобы напомнить это определение и отметить некоторые его преимущества перед <<общепринятым>>.

\section{Преобразования графиков функций}

\begin{exa} Цитируем \cite[С. 53]{cherk}:
\begin{quote}
\textit{Построение графика функции} $y=mf(kx+a)+b$.  [С]начала следует построить график функции $y=f(x)$; затем сжать (или растянуть) его вдоль оси $Ox$ в $k$ раз; перенести влево (или вправо, в зависимости от знака) на $a$; растянуть (или сжать) вдоль оси $Oy$ в $m$ раз и поднять или опустить на~$b$. Если хотя бы один из коэффициентов $k$ или $m$ отрицателен, то следует симметрично отобразить график относительно оси $Oy$ или $Ox$ соответственно.\\
{\sl Замечание.} Описанный порядок действий не является единственным. 
\end{quote}
\end{exa}

Авторы --- преподаватели механико-математического факультета МГУ, широко известные своими многочисленными методическими указаниями и справочными пособиями по элементарной математике для абитуриентов и школьников. И, тем не менее, описанный ими способ построения графика функции ошибочен. Мы говорим об этом не для того, чтобы упрекнуть именитых авторов, а чтобы показать, что даже такие авторитетные специалисты могут ошибаться в элементарных вопросах, если они пренебрегают доказательствами.

Займёмся теперь получением правильного решения --- не столько для того, чтобы исправить ошибку, сколько для того, чтобы показать, как {\em доказательство: а) ведёт к открытию, б) обосновывает его и в) помогает понять полученный результат}.

\begin{exa}[Продолжение примера 1]
Пусть $f$ --- числовая функция, $a$, $b$, $k$, $m$ --- действительные числа и $k\ne 0$, $m\ne 0$. Положим $u(x)= mf(kx+a)+b$. Выясним, как связаны между собой графики~$\Gamma_f$ и~$\Gamma_u$ функций $f$ и $u$.

Чтобы сократить формулировки, для любого $\alpha\in\R$ определим преобразования $H^1_\alpha$, $H^2_\alpha$ плоскости $\R^2$ условиями: $H^1_\alpha(x,\ y)=(\alpha x,\ y)$, соотв., $H^2_\alpha(x,\ y)=(x,\ \alpha y)$, $(x,\ y)\in\R^2$. Для любой точки $(x,\ y)\in\R^2$ имеем:
\begin{align*}
(x,\ y)\in \Gamma_u &\iff y=u(x);\\
                    &\iff y=mf(kx+a)+b;\\
                    &\iff \tfrac{y-b}{m}=f(kx+a);\\
                    &\iff (kx+a,\ \tfrac{y-b}{m})\in\Gamma_f;\\
                    &\iff (kx+a,\ y-b)\in H^2_m\Gamma_f;\\
                    &\iff (kx+a,\ y)\in H^2_m\Gamma_f+(0,\ b);\\
                    &\iff (kx,\ y)\in H^2_m\Gamma_f+(0,\ b)-(a,\ 0);\\
                    &\iff (kx,\ y)\in H^2_m\Gamma_f+(-a,\ b);\\
                    &\iff (x,\ y)\in H^1_\frac1k(H^2_m\Gamma_f+(-a,\ b));\\
                    &\iff (x,\ y)\in H^1_\frac1k H^2_m\Gamma_f+(-\tfrac ak,\ b).
\end{align*}
Таким образом,
\[
  \Gamma_u = H^1_\frac1k H^2_m\Gamma_f+(-\tfrac ak,\ b).
\]
Это означает, что \textit{для получения графика функции
\[
  u\colon x\leadsto y=mf(kx+a)+b
\]
из графика функции $f$ нужно абсциссу каждой точки графика функции~$f$ разделить на $k$, ординату умножить на $m$ и сдвинуть все полученные точки на вектор $(-\tfrac ak,\ b)$}.
\end{exa}

Теперь легко указать ошибку в цитированном пособии. Сжав (или растянув) график исходной функции в $k$ раз вдоль оси $Ox$, авторы пособия переносили его (влево или вправо) на $a$, а нужно было --- на $\frac{a}{k}$. Разумеется, причина ошибки заключается не в недостаточной осведомлённости авторов, а в отсутствии доказательства и установки на него.

Заметим, что и этот пример обладает отмеченными выше необходимыми особенностями: он учит преподаванию математики (не только в отношении построения графиков функций, но и в отношении установки на поиск доказательств), и делается это не в ущерб основному предмету~--- математическому анализу, а с пользой для него.

\section{Работа с кванторами в курсе математического анализа}

Хорошо известно, какие трудности вызывает начало изучения теории пределов первокурсниками. Особенно удручает их нагромождение кванторов в определении предела и доказываемых предложениях: для любого \ldots существует \ldots такое, что для  любого \ldots. А ведь этих трудностей можно было бы полностью или почти полностью избежать при внесении соответствующих изменений в практику преподавания математики в школе, потому что в школьной математике имеется большое число предложений, неявно имеющих такую же кванторную приставку, как и определение предела.

Рассмотрим, например, теорему о медианах треугольника. В традиционной формулировке она говорит, что медианы треугольника пересекаются в одной точке. На вопрос, как понимать слова \glqq в одной точке\grqq, практически все опрашиваемые отвечают: в одной --- значит не в двух, не в трёх\dots. Иначе говоря, они полагают, что слово \glqq одна\grqq\ в этой формулировке является числительным. Тогда мы предлагаем им взглянуть на следующий чертёж и говорим: вот треугольник, вот его медианы, а вот --- одна точка, но медианы почему-то не хотят через неё проходить. И этот незатейливый пример очень многих приводит в замешательство.

При внимательном его рассмотрении выясняется, что дело заключается не в том, чтобы точка была одна (точек пересечения не больше одной даже для двух медиан), а в том, чтобы это была нужная точка. Так же как слово \glqq однажды\grqq\ совсем не обязательно означает \glqq единожды\grqq, а слова \grqq в один прекрасный день\grqq\ отнюдь не означают, что этот день должен быть непременно один, так и здесь слова \glqq в одной точке\grqq\ являются указанием не на единственность точки, а на её {\em существование}: медианы треугольника пересекаются \glqq в одной прекрасной точке\grqq. Это означает, что существует такая точка, которая является общей всем медианам.

Разобравшись с точкой, задумываемся о медианах: какие медианы пересекаются \glqq в одной прекрасной точке\grqq? Та и эта? Какие-нибудь две? Некоторые? Или все? Выясняется, что --- все.

Наконец, переходим к треугольнику: о каком треугольнике речь? О нашем близком знакомом? О каком-то особенном? Нет, обо всех. В итоге этих обсуждений выясняется, что подлинный смысл теоремы о медианах треугольника заключается  в том, что
\begin{center}
  {\bfseries\em для любого} треугольника {\bfseries\em существует} точка,\\ через которую проходят {\bfseries\em все} его медианы.
\end{center}
Без выяснения этого смысла теорема не может быть правильно понята и, тем более, правильно доказана. Значит, все описанные обсуждения служат совершенствованию преподавания математики в школе.

С другой стороны, кванторная приставка у полученного предложения в точности такая же, как и у определения предела, непрерывности, \ldots. Выявление же этого обстоятельства способствует совершенствованию преподавания математического анализа. Особенно, после того как выясняется, что подобных предложений в школьной математике существует очень много, их просто не замечают, потому что не вникают в их логическое строение.

Помимо теоремы о медианах треугольника сюда относятся аналогичные теоремы о биссектрисах треугольника, о серединных перпендикулярах к сторонам треугольника, о высотах треугольника, теоремы об окружностях, вписанной в треугольник и описанной около него, теоремы о сферах, вписанной в тетраэдр и описанной около него, предложение о том, что через каждую прямую в пространстве проходит хотя бы одна плоскость (означающее, что для любой прямой в пространстве существует плоскость, которой принадлежат все точки прямой), и так далее.

Все перечисленные и многие другие предложения имеют одинаковые кванторные приставки и по этой причине одинаковую логическую структуру доказательств: рассмотрим произвольный \ldots, найдём \ldots\ такой, что для всех \ldots. После многократной отработки на знакомом школьном материале, доказательства подобных теорем теряют в глазах студентов пугающую обособленность, становятся в один ряд с привычными школьными утверждениями, благодаря чему студенты обретают необходимую уверенность при работе с подобными теоремами и задачами, как в школьном курсе математики, так и в курсе математического анализа.


\begin{thebibliography}{99}

\bibitem{gned}
{\sc Б.~В.~Гнеденко.} Об образовании преподавателя математики средней школы // Математика в школе. --- 1989, № 3. --- С. 19--22.

\bibitem{what}
{\sc А.~Х.~Назиев.} Чему учить будущего учителя математики? // Известия Российской Академии Естественных Наук. Дифференциальные уравнения. --- 2006, 11. --- С. 170--174.

\bibitem{abma}
{\em Алгебра и начала математического анализа} : учеб. для 10—11 кл. общеобразоват. учреждений / [А.~ Н.~Колмогоров,
А.~М.~Абрамов, Ю.~П.~Дудницын и др.] ; под ред. А.~Н.~Колмогорова. — 17-е изд. — М. : Просвещение,
2008. — 384 с.

\bibitem{func-reviz}
{\sc А.~Х.~Назиев.} Ещё раз о том, что такое функция. // Современные подходы к оценке и качеству математического образования в школе и вузе: Материалы XXXII Международного научного семинара преподавателей математики университетов и педагогических вузов, 26--28 сентября 2013 г. --- Екатеринбург, 2013. --- С. 184--186

\bibitem{humed}
{\sc А.~Х.~Назиев.} {\em Гуманитарно ориентированное преподавание математики в общеобразовательной школе}. Монография // Рязань: Изд-во РИРО, 1999. --- 112 с.

\bibitem{dissB}
{\sc А.~Х.~Назиев.} {\em Гуманитаризация основ специальной подготовки учителей математики в педагогических вузах}: Дисс. \ldots\ докт. пед. наук. --- Москва, МПГУ, 2000. --- 386 с.

\bibitem{Polya:disc}
{\sc Д. Пойа.} {\it Математическое открытие.} --- М.: "<Наука">, 1970.

\bibitem{cherk}
{\sc О.~Ю.~Черкасов, А.~Г.~Якушев.} {\it Математика на вступительных экзаменах в серьёзные вузы.} М.: Московский лицей, 1998.

\bibitem{joung}
{\sc Юнг~В.~А.} {\it Как преподавать математику}.
 --- М.--Пд.: Госиздат, 1923. --- 302~с.

\end{thebibliography}
\end{document}